\documentclass[a4paper,11pt]{article}
\usepackage{hyperref}
\usepackage{mathrsfs}
\usepackage{stmaryrd}
\usepackage{amssymb}
\usepackage{amsfonts}
\usepackage{amscd}
\usepackage{latexsym}
\usepackage{amsmath}
\usepackage{indentfirst}
\begin{document}

\title{\textbf{Approximated by finite-dimensional homomorphisms into Simple C*-Algebras with Tracial Rank One}}
\author{Liu Junping$^1$ \ \ \ Zhang Yifan$^2$  \\[6pt]\\
\footnotesize{1: East China Normal University, Shanghai, China;
jpliu@math.ecnu.edu.cn;}\\
 \footnotesize{2: Xiamen University of
Technology, Xiamen, China; yfzhang@xmut.edu.cn.}}
\date{}
\maketitle

\begin{abstract}

We discuss when a unital homomorphism $\phi: C(X)\rightarrow A$ can be approximated by finite-dimensional homomorphisms,
where $X$ is a compact metric space and $A$ is unital simple C*-algebra with tracial rank one. In this paper, we will give
a necessary and sufficient condition.

\end{abstract}

\section*{1. Introduction}

\noindent In the theory of C*-algebras, studying homomorphisms between two C*-algebras is of fundamental importance. As a simple step,
but also important, we study homomorphisms from some commutative C*-algebra $C(X)$, where $X$ is a compact metric space, into some simple
C*-algebra. Among these homomorphisms, the ones defined by evaluation at some finite points in $X$ are the most simple case, or equivalently,
the ones with finite-dimensional range (we call them finite-dimensional homomorphisms). Now, it is natural to study their limits (in the point-wise convergence topology).

\par In this paper, the target C*-algebra of a homomorphism we shall consider is in an important
class of simple C*-algebras in the classification theory, the unital
simple C*-algebras with tracial rank no more than one. It is introduced by H. Lin to
aid the program of classification of nuclear C*-algebras
(\cite{Lin2}). H. Lin
completely classified the unital nuclear separable simple
C*-algebras with tracial rank one which satisfy the UCT, see
\cite{Lin4}. Let $A$ be a unital simple C*-algebra with tracial rank no more than one, consider a unital monomorphism $\phi: C(X)\rightarrow A$.
(For this problem, we only need to consider monomorphism.)
When $X$ is path-connected and $A$ is of tracial rank zero, it is
proved that $\phi$ can be approximated by finite-dimensional homomorphisms if and only if $\phi$ induces an zero element in
$KL(C_0(X), A)$ (see \cite{HLX}). In the present paper, we shall extend the result to the tracial rank one case.
It is worth mentioning that the latter C*-algebras are not of real rank zero.
It turns out that $\phi$ can be approximated by finite-dimensional homomorphisms if and only if  $[\phi]$ vanishes on $\underline{K}(C_0(X))$,
and in addition, the induced maps $\widehat{\phi}$ maps $\mbox{Aff}T(C(X))$ into $\overline{\rho_A(K_0(A))}$ and $\phi^\dag$ is trivial.
For a general compact metric space $X$ (not necessarily path-connected nor a disjoint union of finitely many path-connected spaces), we need some new
generalized notation
to describe $[\phi]$ (see Definition 3.4).

\par In the literature, this problem is related to the properties such as real rank zero, (FU) and (FN), corresponding to $X=[0,1], \mathbb{T}$, and a compact subset in the complex plane, respectively (\cite{Lin1} and \cite{Lin8}). In \cite{Lin8}, it is shown that a unitary in a unital simple C*-algebra with tracial rank one can be approximated
by unitaries with finite spectrum if and only if $u\in CU(A)$ and $\widehat{u^n+(u^n)^*}, \widehat{i(u^n-(u^n)^*)} \in \overline{\rho_A(K_0(A))}$ for all $n\geq 1$.
As an application of the main result in this paper , we shall describe when a normal element in a unital simple C*-algebra with tracial rank one can be approximated by normal elements with finite spectrum.

\

\section*{2. Preliminaries}

\par In this section, we gather some notations and well-known facts.

\

\noindent \textbf{2.1.} Let $A, B$ be two C*-algebras and let $\phi, \psi: A\rightarrow
B$ be two maps. Suppose that $\mathcal{F}\subset A$ and $\epsilon>
0$. We write
$$\phi\approx_\epsilon\psi \ \mbox{on}\ \mathcal{F}$$
if $\parallel \phi(x)-\psi(x)\parallel<\epsilon$ for all $a\in
\mathcal{F}$. Similarly, we write
$$\phi\approx_\epsilon \mbox{ad}u\circ\psi \ \mbox{on}\ \mathcal{F}$$
if there is a unitary $u\in B$ such that $\parallel
\phi(x)-u\psi(x)u^*\parallel<\epsilon$ for all $a\in \mathcal{F}$.

\

\noindent \textbf{2.2.}  If $X=X_1\sqcup\cdots\sqcup X_m$ is a disjoint union of path-connected compact metric spaces with
each component $X_i$ a
base point $x_i\in X_i$ for $i=1,\cdots,m$. We shall use the notation $C_0(X)$ to mean
the set of continuous functions on $X$ which vanish at all $x_i$.  Put $rC(X)=C(X)/C_0(X)\cong \mathbb{C}^m$ (\cite{DN} and \cite{EG}).

Let $A$ be a unital C*-algebra and $\phi:C(X)\rightarrow A$ be a
unital homomorphism. Then $\phi$ defines an element $[\phi]$ in
$KK(C(X),A)$. It is known that the short exact sequence
$$0\rightarrow C_0(X)\rightarrow C(X) \rightarrow rC(X) \rightarrow 0$$
is split and there is a natural decomposition,
$$KK(C(X),A)=KK(C_0(X),A)\oplus KK(rC(X),A).$$
If $\beta \in KK(C(X),A)$, we will write $\beta=(\beta_0, \beta_1)$
under this decomposition. In particular, suppose that $\phi:
C(X)\rightarrow A$ is a unital homomorphism and denote by $e_i\in C(X)$ the identity of $C(X_i)$, then
$[\phi]=([\phi]_0,[p_1],\cdots,[p_m])\in KK(C_0(X),A)\oplus
KK(rC(X),A)=KK(C_0(X),A)\oplus K_0(A)\oplus\cdots\oplus K_0(A)$, where $p_i=\phi(e_i), i=1,\cdots,m$.

\par From the universal coefficient theorem (see \cite{RS} or
\cite{Bl}), for the C*-algebras $C=C(X)$ and $A$ as above, there is a
split short exact sequence
$$0 \rightarrow \mbox{Ext}^1_{\mathbb{Z}}(K_*(C), K_*(A))\rightarrow KK(C,A) \xrightarrow {\gamma} \mbox{Hom}^0(K_*(C), K_*(A))\rightarrow 0.$$
Define
$$KL(C, A)=KK(C, A)/\mbox{Pext}^1_{\mathbb{Z}}(K_*(C), K_*(A)).$$

\

\noindent \textbf{2.3.} Let $A$ be a C*-algebra and let $C_n$ be a
commutative C*-algebra with $K_0(C_n)=\mathbb{Z}/n$ and
$K_1(C_n)=0$. We use
the following notation:
$$K_*(A,\mathbb{Z}/n)=K_*(A\otimes C_n)$$
and $$\underline{K}(A)=K_*(A)\oplus\bigoplus_{n=2}^\infty
K_*(A,\mathbb{Z}/n).$$ Denote by
$Hom_{\Lambda}(\underline{K}(A),\underline{K}(B))$ the set of
systems of group homomorphisms which is compatible with all the
Bockstein Operations (see \cite{DL} for details). From \cite{DL},
we know that
$$KL(A,B)=Hom_{\Lambda}(\underline{K}(A),\underline{K}(B))$$
for A is a separable amenable C*-algebra which satisfies the UCT.

\par If $\alpha \in KL(C(X),B)$, we use the notation $\alpha=\{\alpha_n^i\}$,
 where $\alpha_n^i\in Hom(K_i(C(X),\mathbb{Z}/n)\rightarrow K_i(B,\mathbb{Z}/n)), i=0,1$.

\

\noindent \textbf{2.4.} Let $X$ be a compact metric space and $x,y\in X$. Let $\delta>0$. We write $x\sim_\delta y$, if there are
points $x_0,x_1,\cdots,x_m$ in $X$ such that
$$x_0=x, x_m=y,\ \mbox{and}\ dist(x_i,x_{i+1})<\delta, $$
for $i=0,1,\cdots,m-1$. A subset $Y\subset X$ is said to be $\delta$-connected, if for any two points $x, y$ in Y, one has $x\sim_\delta y$.
\par It is well-known that for any $\delta>0$, $X$ can be divided into finitely many disjoint $\delta$-connected components. We will
frequently use the following lemma.

\

\noindent \textbf{Lemma} (Lemma 3.3 of \cite{HLX}) \emph{Let $X$ be a compact metric space and $G\subset \underline{K}(C(X))$ be a finitely generated subgroup.
Then there exists $\delta>0$ satisfying the following:}
\par \emph{If $X=X_1\sqcup\cdots\sqcup X_m$, where each $X_i$ is $\delta$-connected, $A$ is a unital C*-algebra, and $\phi, \psi: C(X) \rightarrow A$ are two unital
finite-dimensional homomorphisms such that $[\phi]([e_i])=[\psi]([e_i])$ for $i=1,\cdots,m$,
then} $$[\phi]\mid_G=[\psi]\mid_G.$$ \hspace*{\fill} $\Box$

\

\noindent \textbf{2.5.} For a unital C*-algebra $A$, let $T(A)$
denote the space of all tracial states of $A$. It is well known that
$T(A)$ is a  Choquet simplex. Let $\mbox{Aff}T(A)$ be the space of
all affine continuous real functions on $T(A)$. Then
$\mbox{Aff}T(A)$ is an ordered Banach space with order unit. If $X$
is a compact Hausdorff space, then it is well known that
$\mbox{Aff}T(C(X))=C_{\mathbb{R}}(X)$, the space of all real
continuous functions on $X$.

\

\noindent \textbf{2.6.}  Let $\phi: C\rightarrow A$ be a unital
homomorphism. Denote by $\phi^T: T(A)\rightarrow T(C)$ the affine
continuous map induced by $\phi$, i.e. $\phi^T(\tau)=\tau\circ\phi$
for all $\tau$ in $T(A)$. It then also induces a unital positive
linear map $\widehat{\phi}: \mbox{Aff}T(C)\rightarrow\mbox{Aff}T(A)$
defined by $\widehat{\phi}(f)=f\circ\phi^T$ for all $f$ in
$\mbox{Aff}T(C)$.
\par Let $A$ be a unital C*-algebra, denote by $(K_0(A), K_0(A)^+, [1_A])$ the associated
scaled ordered group and $SK_0(A)$ the state space of $K_0(A)$ (see
\cite{R}). There is an affine continuous map $r_A: T(A)\rightarrow
SK_0(A)$ defined by $r_A(\tau)([p])=\sum_{i=1}^n\tau(p_{ii})$, where
$\tau \in T(A)$ and $[p] \in K_0(A)$ is the element presented by the
projection $p \in M_n(A)$. Then $r_A$ defines a canonical map
$\rho_A: K_0(A)\rightarrow \mbox{Aff}T(A)$ by
$\rho_A(g)(\tau)=r_A(\tau)(g)$ for $\tau$ in $T(A)$ and $g$ in
$K_0(A)$.

\par Let $\pi_A: \mbox{Aff}T(A)\rightarrow \mbox{Aff}T(A)/\overline{\rho_A(K_0(A))}$ denote the canonical quotient map.

\

\noindent \textbf{2.7.} Let $U(A)$ be the unitary group of $A$, and $CU(A)$ the closure of the commutator subgroup
of $U(A)$. Using de la Harpe-Skandalis determinant, by Theorem 3.2 of \cite{Th}, one has the following splitting short exact
sequence:
$$0\rightarrow \mbox{Aff}T(A)/\overline{\rho_A(K_0(A))}  \rightarrow U_\infty(A)/CU_\infty(A)   \rightarrow K_1(A)  \rightarrow 0. \eqno{(2.1)}$$
We then have
$$U_\infty(A)/CU_\infty(A)\cong \mbox{Aff}T(A)/\overline{\rho_A(K_0(A))}\oplus K_1(A). \eqno{(2.2)}$$

For a unital homomorphism $\phi: C\rightarrow A$, it induces a group homomorphism $\phi^\ddag:U_\infty(C)/{CU_\infty(C)}\rightarrow
U_\infty(A)/CU_\infty(A)$. With respective to the decomposition (2.2), we can view $\phi^\ddag$ as a matrix:
$$
\left[
    \begin{array}{cc}
       \alpha_{11} & \alpha_{12}\\
       \alpha_{21} & \alpha_{22}\\
      \end{array}
   \right]
$$
Here $\alpha_{21}$ is automatically zero by exactness, $\alpha_{11}$ is induced by $\widehat{\phi}$, and $\alpha_{22}$ is the map $\phi_{*1}$.
Denote the rest homomorphism $\alpha_{12}$ by
$$\phi^\dag: K_1(C)\rightarrow \mbox{Aff}T(A)/\overline{\rho_A(K_0(A))}.$$

\

\noindent \textbf{2.8.} We shall recall the definition and some
basic properties of unital simple C*-algebras of tracial ranks (see
\cite{Lin2} and \cite{Lin4}).
\par Let $A$ be a unital simple
C*-algebra and $k \in \mathbb{N}$. We say that $A$ has
\textbf{tracial rank no more than $k$} if for any finite subset
$\mathcal{F}\subset A$, any $\epsilon>0$ and any non-zero positive
element $a\in A_+$, there exist a non-zero projection $p\in A$ and a
C*-subalgebra $B$ of $A$ with $1_B=p$ such that
\par (0) $B$ has form $B=\oplus_{i=1}^q P_{n_i}M_{n_i}(C(X_i))P_{n_i}$, where $P_{n_i}$ are projections
in $M_{n_i}(C(X_i))$ and $X_i$ is a finite CW complex with
$dim(X_i)\leq k$ for each $i$,
\par (1) $\parallel px-xp\parallel<\epsilon$ for all $x\in
\mathcal{F}$,
\par (2) $pxp \in_\epsilon B$ for all $x\in
\mathcal{F}$,
\par (3) $[1-p]\leq [a]$.
\par We will write $TR(A)\leq k$ if $A$ has tracial rank no more than $k$.
Especially, $TR(A)=0$ means that $A$ has tracial rank zero, i.e. the
above $B$ can be chosen to be finite dimensional.

\par Recently, H. Lin proved that if a unital simple C*-algebra $A$ with $TR(A)\leq k$ satisfies the UCT, then $A$
actually has tracial rank no more that one (see \cite{Lin11}). Hence we focus on $A$ with $TR(A)\leq1$.
\par Suppose that $A$ is a unital simple C*-algebra with $TR(A)\leq 1$, then it is well known that
$A$ has stable rank one, real rank no more than one, weakly
unperforated $K_0$-group with Riesz interpolation property and
fundamental comparison property (see \cite{Lin4}). For $TR(A)=0$
case, we know that $A$ has real rank zero and the canonical map
$\rho_A: K_0(A)\rightarrow \mbox{Aff}T(A)$ in 2.6 has dense range
when $A$ is infinite dimensional simple C*-algebra.

\

\section*{3. Main Results}

\

\par First of all, we discuss the case when $X$ is path-connected. The following result is the main theorem in \cite{HLX}.

\

\noindent \textbf{Lemma 3.1.}  (Theorem 3.9 of \cite{HLX})\quad
\emph{Let $X$ be a compact path-connected metric space, and $A$ be a unital
simple C*-algebra with $TR(A)=0$. Suppose that $\phi:
C(X)\rightarrow A$ is a unital monomorphism. Then $\phi$ can be approximated by finite-dimensional homomorphisms if and
only if} $$[\phi]_0=0\ \mbox{in}\ KL(C_0(X), A).$$ \hspace*{\fill} $\Box$

\

\par By using the uniqueness theorem recently proved by H. Lin (\cite{Lin10}) and the method used in \cite{Lin8}, we obtain the following result.

\

\noindent \textbf{Theorem 3.2.}  \emph{Let $X$ be a compact path-connected metric space with a based point $x_0$, and $A$ be a unital simple infinite dimensional C*-algebra with $TR(A)\leq 1$.
Suppose $\phi: C(X)\rightarrow A$ is a unital monomorphism. Then $\phi$ can be approximated by finite-dimensional homomorphisms
if and only if} $$[\phi]_0=0\ \mbox{in}\ KL(C_0(X), A),$$
$$ \pi_A\circ \hat{\phi}= 0, \ \mbox{and}  $$
$$\phi^\dag=0.$$

\

\noindent \textbf{Proof.} Firstly, suppose that $\psi: C(X)\rightarrow A$ is a unital finite-dimensional homomorphism. Then we can write $\psi$
as
$$\psi(f)=\sum_{k=1}^mf(x_k)p_k$$
for all $f\in C(X)$, where $x_k\in X$ and $p_1, \ldots, p_m$ are mutually orthogonal projections in $A$ with
$\sum_{k=1}^mp_k=1$. Define $\psi_0: C(X)\rightarrow A$ by
$$\psi_0(f)=f(x_0)\cdot1_A$$
for all $f \in C(X)$. Since $X$ is path-connected, then $\psi$ is homotopic to $\psi_0$, and hence
$$[\psi]_0=[\psi_0]_0=0,\ \mbox{in}\ KL(C_0(X), A). $$ Also, since $\overline{\rho(K_0A)}$ is a $\mathbb{R}$-linear subspace of
 $\mbox{Aff}T(A)$ (see Proposition 3.6 of \cite{Lin8}), we see that $\widehat{\psi}$ maps $C_{\mathbb{R}}(X)$ into $\overline{\rho(K_0A)}$.
Next, if $u\in U(C(X))$, then $\mid u(x_k)\mid=1$, we write $u(x_k)=exp(i\theta_k)$, where $\theta_k\in \mathbb{R}$ for $k=1,\cdots,m$. Put
$h=\sum_{k=1}^m\theta_kp_k\in A_{sa}$, then
$$\psi(u)=\sum_{k=1}^mu(x_k)p_k=exp(ih).$$
Note that $\widehat{h}\in \overline{\rho(K_0A)}$, $\psi(u)\in CU(A)$
by Theorem 2.9 of \cite{Lin8}. Hence $\psi^\ddag(U(C(X)))\subset
CU(A)$. Similarly, $\psi^\ddag(U_n(C(X)))\subset CU_n(A)$. Then
$\psi^\ddag=0$, and hence $\psi^\dag=0.$ Now if $\phi$ can be
approximated by finite-dimensional homomorphisms, $\phi$ must
satisfies the mentioning three conditions.

\par Conversely, suppose that $\phi: C(X)\rightarrow A$ is a unital homomorphism such that $[\phi]_0=0$, $\pi_A \circ \widehat{\phi}=0$,
and $\phi^\dag=0$. As the proof of Lemma 4.1 of \cite{Lin8}, we can choose a unital simple C*-subalgebra $B\subset A$ with tracial rank zero such
that the inclusion $\imath: B\rightarrow A$ induces an isomorphism:
$$(K_0(B), K_0(B)_+, [1_B], K_1(B)) \cong (K_0(A), K_0(A)_+, [1_A], K_1(A)).$$
Then $[\imath]$ is a KK-equivalence in $KK(B, A)$ (see 7.6 of \cite{RS}), and hence there is a $\beta\in KK(A,B)$ such that
$$ [\imath]\times \beta=[id_B]\ \mbox{and}\ \beta\times[\imath]=[id_A].$$
Define $\kappa=[\phi]\times \beta\in KL_e(C(X), B)^{++}$. Since $TR(B)=0$, we know that $\mbox{Aff}T(B)=\overline{\rho_B(K_0(B))}=\overline{\rho_A(K_0(A))}$.
By assumption, $\widehat{\phi}$ maps $C_{\mathbb{R}}(X)$ into $\overline{\rho_A(K_0(A))}$, then $\widehat{\phi}$ defines a unital strickly
positive linear map $\gamma: C_{\mathbb{R}}(X)\rightarrow \mbox{Aff}T(B)$ such that $\widehat{\imath}\circ\gamma=\widehat{\phi}$. We now check that the
defined pair $(\kappa, \gamma)$ is compatible, that is, $\rho_B\circ\kappa_0^0=\gamma\circ\rho_{C(X)}$. It suffices to show that
$\widehat{\imath}\circ\rho_B\circ\kappa_0^0=\widehat{\imath}\circ\gamma\circ\rho_{C(X)}$, since $\widehat{\imath}$ is injective. This is equivalent to
$$\rho_A\circ\widehat{\imath}\circ\kappa_0^0=\widehat{\phi}\circ\rho_{C(X)}.\eqno{(3.1)}$$
Since $\widehat{\imath}\circ\kappa_0^0=\phi_{*0}$, the equation (3.1) becomes
$$\rho_A\circ\phi_{*0}=\widehat{\phi}\circ\rho_{C(X)},$$
and this is well-known true. The following diagram shows the above calculation:

\[
\begin{CD}
K_0(C(X))  @>\kappa_0^0>> K_0(B)@> \xleftarrow{\beta_0^0}>\imath_{*0}>K_0(A)\\
 @VV{\rho_{C(X)}}V @VV{\rho_{B}}V @VV{\rho_{A}}V\\
C_{\mathbb{R}}(X) @>>\gamma> \mbox{Aff}T(B))@>>\widehat{\imath}>\mbox{Aff}T(A)\\
\end{CD}
\]

Now apply Theorem 5.2 of \cite{Lin6}, there exists a unital monomorphism $h:C(X)\rightarrow B$  such that
$$[h]=\kappa\ \mbox{in}\ KL(C(X), B)\quad \mbox{and}\ \widehat{h}=\gamma.\eqno{(3.2)}$$
Then
$$[\imath\circ h]=[h]\times[\imath]=\kappa\times[\imath]=[\phi]\times\beta\times[\imath]=[\phi],\eqno{(3.3)}$$
and $$\widehat{(\imath\circ h)}=\widehat{\imath}\circ \widehat{h}=\widehat{\imath}\circ \gamma=\widehat{\phi}.\eqno{(3.4)}$$
Note that $TR(B)=0$, $h^\dag$ is automatically zero and then
$$(\imath\circ h)^\dag=\phi^\dag=0. \eqno{(3.5)}$$
Combine $(3.3)-(3.5)$ and use Theorem 5.10 of \cite{Lin10}, we conclude that the two homomorphisms $\phi$ and $\imath\circ h$ are approximately unitarily equivalent.
Finally, since $[\phi]_0=0$ in $KL(C_0(X), A)$, it is obviously that $[h]_0=0$ in $KL(C_0(X)), B)$. From Lemma 3.1, $h$ can be approximated by
finite-dimensional homomorphisms, and so is $\phi$.
\hspace*{\fill}$\Box$

\

\noindent \textbf{Corollary 3.3.} \emph{Let $X$ be a disjoint union of finitely many compact
 path-connected metric spaces, and $A$ be a unital simple infinite dimensional
 C*-algebra with $TR(A)\leq 1$.
Suppose $\phi: C(X)\rightarrow A$ is a unital monomorphism. Then $\phi$ can be approximated by finite-dimensional homomorphisms
if and only if} $$[\phi]_0=0\ \mbox{in}\ KL(C_0(X), A),$$
$$ \pi_A\circ \hat{\phi}= 0, \ \mbox{and}  $$
$$\phi^\dag=0.$$

\

\noindent \textbf{Proof.} Write $\phi$ as a finite direct sum and apply Theorem 3.2.\hspace*{\fill}$\Box$

\

\par Now we turn to consider a general (not necessarily path-connected) compact metric space $X$. We need a condition which
generalize the first condition of Corollary 3.3. For any $x\in X$, the evaluation map $\pi_x: C(X)\rightarrow \mathbb{C}$
defines an element $[\pi_x]: \underline{K}(C(X))\rightarrow \underline{K}(\mathbb{C})$.

Let $A$ be a unital C*-algebra and $\phi: C(X)\rightarrow A$
be a unital homomorphism. We shall use the following equation
$$[\phi]\mid_{\cap_{x\in X}Ker[\pi_x]}=0,\eqno{(3.6)}$$
to describe $\phi$. It is easy to see that (3.6) is equivalent to say $[\phi]_0=0$ when $X$ is a disjoint union of finitely
many path-connected compact spaces.

\par Let $X$ be a compact metric space. It is well-known that there is an inductive system $\{C(X_n), h_n\}$ such that $C(X)=\underrightarrow{\lim}{(C(X_n), h_n)}$, where each $X_n$ is a finite CW complex. Suppose that $h_{n, \infty}: C(X_n)\rightarrow C(X)$ is the induced homomorphism. Also there are
 continuous maps $f_n: X\rightarrow X_n$ induce $h_{n,\infty}$. Write $X_n=X_n^1\sqcup\cdots\sqcup X_n^{r(n)}$, where $X_n^i$ is the path-connected component of $X_n$ for
$i=1,\cdots,r(n)$. Fix a point $x_n^i$ in each $X_n^i$.

\

\noindent \textbf{Lemma 3.4.} \emph{With above notations, let $A$ be a unital C*-algebra and $\phi: C(X)\rightarrow A$ be a unital homomorphism.
If (3.6) holds, then we have}
$$[\phi\circ h_{n,\infty}]\mid_{\cap_{i=1}^{r(n)}ker[\pi_{x_n^i}]}=0,\eqno{(3.7)}$$
\emph{for all $n$.}

\

\noindent \textbf{Proof.} Fix a number $n$. For any $g\in\cap_{i=1}^{r(n)}ker[\pi_{x_n^i}]$, let $G$ be the subgroup of $\underline{K}(C(X))$ generated
by $[h_{n,\infty}](g)$.  Let $\delta>0$ be as required by Lemma 2.4 for $G$.  Write $X=Y_1\sqcup\cdots\sqcup Y_N$, where each $Y_k$ is a $\delta$-connected component of $X$.
For each $1\leq k\leq N$, there is $1\leq i(n,k)\leq r(n)$ such that $f_n(Y_k)\cap X_n^{i(n,k)}\neq \emptyset$.
Choose points $y_k \in Y_k$ and $z_n^{i(k)} \in X_n^{i(n,k)}$ such that
$$f_n(y_k)=z_n^{i(k)},$$
for each k. Then
$$[\pi_{y_k}]([h_{n,\infty}](g))=[\pi_{z_n^{i(k)}}](g)=[\pi_{x_n^{i(k)}}](g)=0.\eqno{(3.8)}$$
By Lemma 2.4, we obtain that $[\pi_y]([h_{n,\infty}](g))=0$ for all $y\in Y_k$. This is true for each $k$ and then for all $y\in X$.
Hence, (3.7) is now followed from (3.6). \hspace*{\fill}$\Box$

\

\par The key point of using equation (3.6) is shown in the next lemma.

\

\noindent \textbf{Lemma 3.5.}\emph{ Let $\phi:C(X)\rightarrow A$ be a unital homomorphism as above, then the following two conditions are equivalent:}

\par (i) $[\phi]\mid_{\cap_{x\in X}Ker[\pi_x]}=0;$

\par (ii)\ \emph{for any finitely generated subgroup $G\subset \underline{K}(C(X))$, there is a unital finite dimensional homomorphism
$\psi: C(X)\rightarrow A$ such that}
$$[\phi]\mid_G=[\psi]\mid_G.$$

\

\noindent \textbf{Proof.} Firstly, we prove ``(i) $\Rightarrow$ (ii)''.
Suppose that $C(X)=\underrightarrow{\lim}{(C(X_n), h_n)}$,
where each $X_n$ is a finite CW complex. Write $X_n=X_n^1\sqcup\cdots\sqcup X_n^{r(n)}$, where $X_n^i$ is the path-connected component of $X_n$ for
$i=1,\cdots,r(n)$. And denote by $e_n^i$ the unit of the summand $C(X_n^i)$.
Then, we can further assume that $h_n$ is unital and $h_n(e_n^i)\neq 0$ for all $n, i$, since otherwise, we can delete the summand $C(X_n^i)$
from the inductive system without changing the inductive limit.
Suppose that $h_{n, \infty}: C(X_n)\rightarrow C(X)$ is the induced map. It follows that $$h_{n,\infty}(e_n^i)\neq 0 \eqno{(3.9)}$$
for all $n, i$.
On the other hand, there are continuous maps $f_n: X\rightarrow X_n$ induce $h_{n,\infty}$. Then (3.9) is equivalent to
$$f_n(X)\cap X_n^i\neq \emptyset \eqno{(3.10)}$$
for $i=1,\cdots, r(n)$.
\par Given a finitely generated subgroup $G\leq \underline{K}(C(X))$, we can choose $n$ large enough such that
$G\subset[h_{n,\infty}](\underline{K}(C(X_n)))$. From (3.10), we can
choose $x_n^i\in X_n^i$ and $y^i\in X$ such that
$$f_n(y^i)=x_n^i \eqno{(3.11)}$$
for each $i=1,\cdots,r(n)$.
Now, consider $\phi\circ h_{n,\infty}:C(X_n)\rightarrow A$.
By Lemma 3.4, we have
$$[\phi\circ h_{n,\infty}]\mid_{\cap_{i=1}^{r(n)}ker[\pi_{x_n^i}]}=0.\eqno{(3.12)}$$
\par Define $\psi: C(X)\rightarrow A$ by
$$\psi(f)=\sum_{i=1}^{r(n)}f(y^i)p_i,$$
where $p_i=\phi\circ h_{n,\infty}(e_n^i)$. From (3.12) and the fact that the homomorphisms are unital, we
obtain that $[\phi\circ h_{n,\infty}]=[\psi\circ h_{n,\infty}]$, and hence $[\phi]\mid_G=[\psi]\mid_G.$

\par Conversely, suppose that the condition (ii) is true. Let $g$ be an element in $\cap_{x\in X}ker[\pi_x]$. Let $G$ be the group generated by
$g$. Now we can choose a unital finite dimensional homomorphism $\psi:C(X)\rightarrow A$ such that
$$[\phi]\mid_G=[\psi]_G.\eqno(3.13)$$
Choose a large $n$ satisfying the conditions as above. Then there is $\widetilde{g}\in \underline{K}(C(X_n))$ such that $g=[h_{n,\infty}](\widetilde{g})$.
Since $g\in \cap_{x\in X}Ker[\pi_x]$, by (3.11) , we have $\widetilde{g}\in \cap_{i=1}^{r(n)}ker[\pi_{x_n^i}]$. Therefore by (3.13),
$$[\phi](g)=[\phi]\circ[h_{n,\infty}](\widetilde{g})=[\psi]\circ[h_{n,\infty}](\widetilde{g})=0.$$
\hspace*{\fill}$\Box$

\

\noindent \textbf{Lemma 3.6.}  \emph{Let $A$ be a unital C*-algebra satisfying the fundamental comparison property,
$\phi:C(X)\rightarrow A$ be a unital homomorphism. If $\widehat{\phi}(C_{\mathbb{R}}(X))\subset \overline{\rho_A(K_0(A))}$,
then for any finite subset $\mathcal{F}\subset C_{\mathbb{R}}(X)$ and $\epsilon>0$, there exists a unital finite dimensional
homomorphism $\psi: C(X)\rightarrow A$ such that}
$$\widehat{\phi}\approx_{\epsilon} \widehat{\psi}\ \mbox{on}\ \mathcal{F}.$$

\

\noindent \textbf{Proof.} The proof is contained in the proof of Theorem 3.6 of \cite{HLX}. \hspace*{\fill}$\Box$

\

\par For our purpose, we need further the following approximated version of the uniqueness theorem in \cite{Lin10}.

\

\noindent \textbf{Lemma 3.7.} (Theorem 5.3 of \cite{Lin10})  L\emph{et $X$ be a compact
metric spaces, and $A$ be a unital simple C*-algebra with $TR(A)\leq 1$.
Suppose $\phi: C(X)\rightarrow A$ is a unital monomorphism. For any finite subset $\mathcal{F}\subset C(X)$ and $\epsilon>0$,
there exist a finitely generated subgroup $G\subset \underline{K}(C(X))$, a finite subset $\mathcal{G}\subset C_{\mathbb{R}}(X)$
and $\gamma>0$ such that}
\par \emph{For any unital homomorphism $\psi: C(X)\rightarrow A$, if}
$$[\phi]\mid_G=[\psi]_G,$$
$$\widehat{\phi}\approx_{\gamma} \widehat{\psi}\ \mbox{on}\ \mathcal{G},\ $$
\emph{and}
$$\phi^\ddag=\psi^\ddag,$$
\emph{then there is a unitary $u\in A$ such that}
\vspace*{2mm}
\par \hspace*{\fill}$\phi\approx _ \epsilon\mbox{ad}u\circ \psi\quad on\quad
\mathcal{F}.$ \hspace*{\fill}$\Box$

\

\par With the above preparations, now we can prove the main lemma.

\

\noindent \textbf{Lemma 3.8.} \emph{Let $X$ be a compact
metric spaces, and $A$ be a unital simple infinite dimensional
C*-algebra with $TR(A)\leq 1$.
Suppose $\phi: C(X)\rightarrow A$ is a unital monomorphism satisfying the following conditions:}

\par (i) \emph{for any finitely generated subgroup $G\subset \underline{K}(C(X))$, there is a unital finite dimensional homomorphism
$\psi: C(X)\rightarrow A$ such that}
$$[\phi]\mid_G=[\psi]\mid_G,$$
\par (ii)$$\phi^\ddag=0.$$
\emph{Then $\phi$ can be approximated by finite-dimensional homomorphisms.}

\

\noindent \textbf{Proof.} It suffices to prove the following: for any finite subset $\mathcal{F}\subset C(X)$ and $\epsilon>0$,
there is a unital finite dimensional homomorphism $\psi: C(X)\rightarrow A$ and a unitary $u\in A$ such that
$$\phi\approx _ \epsilon\mbox{ad}u\circ \psi\quad on\quad
\mathcal{F}.$$
Let $G\subset \underline{K}(C(X))$ be a finitely generated subgroup, $\mathcal{G}\subset C_{\mathbb{R}}(X)$ a finite subset and $\gamma>0$ with the properties
described in Lemma 3.7 corresponding to $\mathcal{F}$ and $\epsilon$. For this finitely generated subgroup $G$, we can choose $\delta>0$ as in Lemma 2.4.
Write $X=X_1\sqcup\cdots\sqcup X_m$, where each $X_k$ is a $\delta$-connected component of $X$.

Now let $G_0$ be the group generated by $G$ and $\{[e_0], \cdots, [e_m]\}$. From the condition (i) we can choose a unital finite
dimensional homomorphism $\psi_0: C(X)\rightarrow A$ such that
$$[\phi]\mid_{G_0}=[\psi_0]\mid_{G_0}. \eqno{(3.14)}$$
Under the identification $C(X)=C(X_1)\oplus\cdots \oplus C(X_m)$, we consider the homomorphism $\phi_i:C(X_i)\rightarrow p_iAp_i$ induced
by $\phi$, where $p_i=\phi(e_i)$ for $i=1,\cdots, m$. Also view $\mathcal{G}=\mathcal{G}_1\oplus\cdots\oplus \mathcal{G}_m$.
Since $\phi^\ddag=0$, we have
$$\hat{\phi}(C_{\mathbb{R}}(X))\subset \overline{\rho_A(K_0(A))}.\eqno{(3.15)}$$
And then
we also have $\hat{\phi_i}(C_{\mathbb{R}}(X_i))\subset \overline{\rho_{p_iAp_i}(K_0(p_iAp_i))}$.
For $\mathcal{G}_i$ and $\gamma>0$, by applying Lemma 3.6, there is a unital finite-dimensional homomorphism $\phi_i^{'}:C(X_i)\rightarrow p_iAp_i,$
such that
$$\widehat{\phi_i}\approx_{\gamma}\widehat{\phi_i^{'}}\ \mbox{on}\ \mathcal{G}. \eqno{(3.16)}$$
Define $\psi: C(X)\rightarrow A$ by
$$\psi=\phi_1^{'}\oplus\cdots\oplus\phi_m^{'}.$$
It follows from Lemma 2.4, (3.14), and the definition of $\psi$ that
$$[\psi]\mid_G=[\psi_0]\mid_G=[\phi]\mid_G.\eqno{(3.17)}$$
From (3.16), we have $$\widehat{\phi}\approx_{\gamma}\widehat{\psi}\ \mbox{on}\ \mathcal{G}.\eqno{(3.18)}$$
Note that $\phi^\ddag=\psi^\ddag=0$. We complete the proof by applying Lemma 3.7. \hspace*{\fill}$\Box$

\

\noindent \textbf{Theorem 3.9.} \emph{Let $X$ be compact
 metric spaces, and $A$ be a unital simple infinite dimensional
 C*-algebra with $TR(A)\leq 1$.
Suppose $\phi: C(X)\rightarrow A$ is a unital monomorphism. Then $\phi$ can be approximated by finite-dimensional homomorphisms
if and only if} $$[\phi]\mid_{\cap_{x\in X}Ker[\pi_x]}=0,$$
$$ \pi_A\circ \hat{\phi}= 0, \ \mbox{and}  $$
$$\phi^\dag=0.$$

\

\noindent \textbf{Proof.} For the necessity, we need only to check the first condition. This is easy by Lemma 3.5. Conversely, combine
Lemma 3.5 and 3.8, and also note that the three conditions imply that $\phi^\ddag=0$, see 2.7. \hspace*{\fill}$\Box$

\

\par As an application, we now obtain a necessary and sufficient condition for a normal element being approximated by normal elements with finite spectrum.

\

\noindent \textbf{Corollary 3.10.} \emph{Let $A$ be a unital simple infinite dimensional C*-algebra with $TR(A)\leq 1$ and
$x\in A$ be a normal element. Suppose that $\phi:C(sp(x))\rightarrow A$ is the map induced by continuous functional calculus.
Then $x$ can be approximated by normal elements with finite spectrum if and only if}
$$\phi^\ddag=0.$$

\

\noindent \textbf{Proof.} Put $X=sp(x)$. We know that $K_*(C(X))$ is torsion free, then
$$KL(C(X), A)=\bigoplus_{i=0,1}Hom(K_i(C(X)), K_i(A)).$$
Note that in this case, $[\phi]_0=0$ is equivalent to say that $\phi_{*1}=0$. Then the three conditions in Theorem 3.9
follow by $\phi^\ddag=0$. \hspace*{\fill}$\Box$

\

\noindent \textbf{Acknowledgement.} \ The authors would like to thank Professor Huaxin Lin for many useful conversations.

\

\def\refname{4. References}

\end{document}